\documentclass[12pt]{amsart}
\usepackage{amsmath,amssymb, amsthm, latexsym, amsfonts}
\usepackage{graphicx}

\newcommand\Z{\mathbb{Z}}

\newtheorem{thm}{Theorem}[section]

\newtheorem{lemma}{Lemma}[section]

\usepackage{multicol}
\usepackage{vmargin, marginnote}
\usepackage[mathscr]{eucal}
\usepackage{xcolor}

\newcommand{\beq}{\begin{equation*}}

\newcommand{\eeq}{\end{equation*}}

\newcommand{\Beq}{\begin{equation}}

\newcommand{\Eeq}{\end{equation}}

\title[Existence and numerical approximation of a Schr\"odinger equation]{Existence and numerical approximation of solutions of a Schr\"odinger equation with derivative in the nonlinear term}
\author[J.C. Mu\~noz Grajales]{Juan Carlos Mu\~noz Grajales }
\author[D. Pizo]{Deissy Marcela Pizo}

\thanks{Juan Carlos Mu\~noz Grajales, email: juan.munoz@correounivalle.edu.co.
Corresponding author.  Deissy Marcela Pizo, email: deissy.pizo@correounivalle.edu.co.
This work was supported by Universidad del Valle under the research project 71299 and MinCiencias under project FP44842-266-2017.  }

\begin{document}

\maketitle

{\footnotesize
\centerline{ Departamento de Matem\'aticas }
\centerline{ Universidad del Valle }
\centerline{Calle 13 Nro 100-00}
\centerline{ Cali -- Colombia }
}

\bigskip

\begin{abstract}
In this paper, we study a Schr\"odinger-type equation featuring a derivative in the nonlinear term and incorporating diffusion effects. This type of equation arises in various physical applications, such as modeling low-order magnetization in ferromagnetic nanocables and describing the collision of ferromagnetic solitons in weakly ferromagnetic media.
We establish a local well-posedness result for the Cauchy problem associated with this model and analyze the convergence and error order of a Fourier spectral numerical scheme for approximating its solutions in the periodic setting. Additionally, we investigate the behavior of solutions in certain asymptotic regimes, both analytically and through numerical experiments, by examining limiting cases of the model's parameters.
\end{abstract}

\section{Introduction}

In this paper, we consider the Schr\"odinger-type equation
\begin{align} \label{schro1}
i \partial_t u + \partial_x^2 u  + i \alpha \partial_x( |u|^2 u ) =  
i ( \eta \partial_x^2 u + \beta \partial_x^3 u + \gamma \partial_x u ), ~~ x \in \mathbb{R}, t\geq 0,
\end{align}
subject to the initial condition 
\begin{equation}\label{initial_cond}
u(x,0)=u_0(x), 
\end{equation}
and $L$-periodic spatial boundary conditions
\begin{equation}\label{periodic_cond}
u(x + L, t) = u(x, t),
\end{equation}
for $x \in \mathbb{R}$, $t \geq 0$. Here $i$ is the complex imaginary unit, $\alpha, \beta, \gamma$ are real constants, $\eta >0$, and $u: \mathbb{R} \times [0, \infty) \to \mathbb{C}$ is a complex function. This equation is a generalization of a mathematical model considered in \cite{Saravanan1}, \cite{Saravanan2}, to study low order magnetization of ferromagnetic nanocables and and collision of ferromagnetic solitons propagating in a weakly ferromagnetic medium, and in \cite{Rogister}, \cite{Bacelli}, \cite{Mj}, \cite{Krause}, as a model to describe the propagation of Alfven waves in plasma physics.

The Schr\"odinger-type equation \eqref{schro1}, in the case where $\eta = \beta = \gamma = 0$, plays a significant role in spatial plasma physics \cite{Mio}, \cite{Ver1}, \cite{Ver2}, and has been used to study nonlinear processes in the Earth's magnetosphere \cite{Pok}. It was first derived in \cite{Mio}, \cite{Mj} to describe the propagation of circularly polarized Alfv\'en waves in plasma. Consequently, this model has found substantial applications in physics. For this reason, the equation has attracted the attention of mathematicians, who have studied it from a theoretical perspective \cite{Tsutsumi1}, \cite{Tsutsumi2}, \cite{Takaoka}, \cite{Kalisch}. In the case where $\eta = \gamma = 0$ and $\beta \neq 0$, equation \eqref{schro1} models the magnetization of a ferromagnetic nanocable, as investigated in \cite{Saravanan1}, \cite{Saravanan2}. On the other hand, when $\beta = 0$, $\eta \neq 0$, and $\gamma \neq 0$, the equation has been studied as a model for Alfv\'en wave propagation \cite{Rogister}, \cite{Bacelli}, \cite{Mj}, \cite{Krause}.

The generalized equation \eqref{schro1} encompasses these physical models, raising important questions regarding the role of the diffusive term (controlled by the parameter $\eta$) and the dispersive term (associated with $\beta$). One of the key problems we address in this work is the theoretical and numerical analysis of how these terms affect the behavior of solutions, particularly in the asymptotic limits $\eta \to 0$ and $\beta \to 0$. To the best of our knowledge, no comprehensive theoretical study exists for equation \eqref{schro1} with arbitrary parameter values. Addressing this gap constitutes the first main contribution of this work.

From a computational perspective, various numerical schemes have been developed for related problems. Finite difference methods for the classical cubic Schr\"o\-din\-ger equation are discussed in \cite{Becerril}, while spectral schemes have been explored in \cite{Kosloff} and more recently in \cite{MunozVargas}. A finite difference scheme for a Schr\"odinger equation with a nonlinear derivative term of the form $|u|^2 \frac{\partial u}{\partial x}$ is proposed in \cite{Cun}, though this differs from the structure of equation \eqref{schro1}.

Numerical methods are essential tools for exploring the properties of solutions to partial differential equations like \eqref{schro1}, particularly in the absence of explicit solutions for general initial data. Unlike the classical Schr\"odinger equation -- whose numerical analysis is well established in the literature -- equation \eqref{schro1} has received comparatively little numerical attention, even in the simplified case $\eta = \beta = \gamma = 0$. We believe that the numerical analysis presented here will contribute to a deeper understanding of the solution space of equation \eqref{schro1}. This forms the second main contribution of the present work.

This paper is organized as follows. In Section 2, we discuss the linear semigroup associated with equation \eqref{schro1} and study the local existence and uniqueness of solutions to the problem \eqref{schro1}-\eqref{periodic_cond}. Section 3 is devoted to the theoretical analysis of the model in the limiting regimes $\eta \to 0$ and $\gamma \to 0$. In Section 4, we introduce a spectral numerical method for approximating the solutions, analyze the convergence error of the semi-discrete formulation, and present numerical experiments to investigate the solution behavior in various parameter regimes.

\section{Existence and uniqueness}
In this section, we study the local well posedness of the nonlinear problem \eqref{schro1}-\eqref{periodic_cond} in an appropriate Sobolev space by using the classical fixed point principle. 

In first place, observe that the problem \eqref{schro1}-\eqref{periodic_cond} can be written as
\begin{equation}\label{schrononlinear}
\begin{aligned}
&\partial_t u(t) = \mathcal{A}(u) + F(u),\\
&u(0) = u_0,
\end{aligned}
\end{equation}
where
\[
\mathcal{A}: = i \partial_x^2 + \eta \partial_x^2 + \beta \partial_x^3 + \gamma \partial_x,
\]
and
\[
F(u) = - \alpha \partial_x ( |u|^2 u ).
\]

\subsection{Linear problem}
In this part, we consider the linear case of problem \eqref{schro1}-\eqref{periodic_cond}
\begin{equation}
\begin{aligned} \label{schro1linear}
&\partial_t u(t) = \mathcal{A}(u), \\
&u(0) = u_0.
\end{aligned}
\end{equation}

In order to construct solutions in this linear case, we consider the Fourier expansion
\begin{equation}
u(x,t) = \sum_k \hat{u}_k(t) e^{i w_k x},
\end{equation}
where $\hat{u}_k$ denotes the $\hat{u}_k, k \in \Z$ denotes the $k$-th Fourier coefficient
\[
\hat{u}_k :=  \frac{1}{L} \int_0^L u(x) e^{-i w_k x}, ~~~ w_k = \frac{2\pi k}{L}.
\]
By substituting this expression into equation \eqref{schro1linear}, we obtain
\begin{equation}
\begin{aligned}
&\hat{u}_k'(t) = (-i w_k^2 - \eta w_k^2 - \beta i w_k^3 + i \gamma w_k) \hat{u}_k(t),\\
&\hat{u}_k(0) = \hat{u}_{0,k}(0).
\end{aligned}
\end{equation}
Therefore, we can solve explicitly for $\hat{u}_k(t)$ to obtain that
\begin{equation*}
\hat{u}_k(t) = \hat{u}_{0,k}\text{exp}\Big( ( -i w_k^2 - \eta w_k^2 - \beta i w_k^3 + i \gamma w_k )t \Big),
\end{equation*}
and thus, we can express the exact solution of problem \eqref{schro1linear} in the form
\[
u(x,t) = \mathcal{T}(t) u_0(x),
\]
where $\mathcal{T}(t), t \geq 0$ is the family of linear operators defined as
\[
\mathcal{T}(t)u(x) := \sum_k  \Big( \hat{u}_{k}\text{exp}( ( -i w_k^2 - \eta w_k^2 - \beta i w_k^3 + i \gamma w_k )t ) \Big) \text{exp} (i w_k x ), ~~ t \geq 0.
\]

\begin{thm}
For $\eta>0$, the family $(\mathcal{T}(t))_{t \geq 0}$ is a $C_0$-semigroup of linear operators and bounded operators in the Sobolev space
\[
H_{per}^s = H_{per}^s([0, L]) := \Big\{ f \in L^2([0,L]): \| f\|_{H^s}^2 = \sum_k (1+ k^2)^s |\hat{f}_k |^2 < \infty \Big\},  
\]
for $s \geq 0$.
For simplicity, we will denote $\| f \|_{s} = \| f \|_{H^s}$. The set $H_{per}^s$
is also a Hilbert space with the inner product
\[
\langle f, g \rangle_s = L  \sum_k (1 + k^2 )^s \hat{f}_k \overline{ \hat{g}_k}.
\]
In particular, when $s = 0$, we get the Hilbert space denoted by
$L_{per}^2 = L_{per}^2(0,L) = H_{per}^0$ with norm given by
\[
\| f \| =\Big(  \int_0^L |f(x)|^2 dx \Big)^{1/2}.
\]
It is important to note that the space $L_{per}^2(0,L)$ is isometrically isomorphic to $L^2 = L^2(0,L)$.

Furthermore, $\mathcal{A}: H_{per}^s \to L^2$ is the infinitesimal generator of the semigroup $(\mathcal{T}(t) ), t \geq 0$,
i.e.
\[
\lim_{t \to 0} \frac{ \mathcal{T}(t) u - u}{t} = \mathcal{A}(u).
\]

\end{thm}
{\bf Proof:} Directly, using the properties of the Fourier transform, it can be shown the following properties of the family $(T(t)), t \geq 0$:

\begin{enumerate}

\item $\mathcal{T}(0) = I$, where $I$ is the identity operator.

\item $\| \mathcal{T}(t) u \|_s \leq C \| u \|_s$, for all $u \in H_{per}^s$, and some positive constant $C$.

\item $\mathcal{T}(t + s) u = \mathcal{T}(t) \mathcal{T}(s) u$, for all $s, t \neq 0$, $u \in H_{per}^s$.

\item The linear operator $\mathcal{A}$ is well defined from $H_{per}^s$ to $L^2$.

\end{enumerate}

\begin{thm}
Let $\lambda>0, \eta >0$. For $u \in H_{per}^s$, $s\geq 0$ we have
\begin{equation}\label{semigroup_prop}
\| \mathcal{T}(t) u \|_{s +\lambda } \leq C_\lambda \sqrt{  1+ \frac{1}{(t \theta )^\lambda} } ~ \| u\|_{s},
\end{equation}
where $\theta = 8 \eta \pi^2/L^2$. Furthermore,
\[
\| \mathcal{T}(t) \partial_x u \|_s \leq C \sqrt{ 1 + \frac{1}{t \theta} } ~ \| u\|_s.
\]

\end{thm}
{\bf Proof:}
Let $u \in H_{per}^s$. Then
\begin{align*}
&\| \mathcal{T}(t) u \|_{s+\lambda}^2 =  \sum_k (1+ k^2 )^{s + \lambda} \Big| e^{ -[ (i+\eta) w_k^2 + \beta i w_k^3 - i \gamma w_k] t }\hat{u}_k \Big|^2 \\
& \leq \sum_k (1+ k^2 )^{s+\lambda} e^{-t \theta k^2} |\hat{u}_k|^2\\
& \leq C \sum_k (1 + k^2 )^s (1 + k^{2 \lambda} ) e^{-t \theta k^2} |\hat{u}_k|^2\\
& \leq \max_{k \in \Z} ( ( 1 + k^{2 \lambda} ) e^{-t \theta k^2} ) \| u\|_s^2\\
& \leq C  \max_{k \in \Z}  ( 1 + k^{2\lambda} e^{-t \theta k^2} ) \| u\|_s^2.
\end{align*}
On the other hand, 
$$
k^{2\lambda} e^{-t \theta k^2} \leq \Big( \frac{\lambda}{t \theta} \Big)^\lambda e^{-\lambda}. 
$$
Consequently,
\[
\| \mathcal{T}(t) u \|_{s + \lambda }^2 \leq C_\lambda \sqrt{  1+ \frac{1}{(t \theta )^\lambda} }~  \| u\|_{s}.
\]
Furthermore, since $\partial_x \mathcal{T}(t) u = \mathcal{T}(t) \partial_x u$, then 
\[
\| \mathcal{T}(t) \partial_x u \|_s = \| \partial_x \mathcal{T}(t) u \|_s \leq \| \mathcal{T}(t) u \|_{s+1}.
\]
Using inequality \eqref{semigroup_prop} with $\lambda = 1$, we obtain
\[
\| \mathcal{T}(t) \partial_x u \|_s \leq C  \sqrt{ 1 + \frac{1}{t \theta} } ~ \| u\|_s.
\]

\subsection{Nonlinear problem}

To establish the existence and uniqueness of solutions of full nonlinear problem \eqref{schrononlinear}, observe that it can be written as
\begin{equation}\label{schrointegral}
u(x,t) = \mathcal{T}(t) u_0(x) + \int_0^t \mathcal{T}(t- \xi) F(u(\xi)) d\xi.
\end{equation}
Thus, if we define the nonlinear integral operator
\[
\Psi(u) := \mathcal{T}(t) u_0(x) + \int_0^t \mathcal{T}(t- \xi) F(u(\xi) ) d\xi,
\]
the problem \eqref{schrointegral} is equivalent to that of finding a fixed point of the
operator $\Psi$.

In first place, we need to analyze the nonlinear term $F(u) = -\alpha \partial_x ( |u|^2 u )$.

\begin{thm}
Let $u, v \in H_{per}^s$ with $s > 1/2$. Then
\begin{equation}
\| |u|^2 u - |v|^2 v \|_{s} \leq C( \| u\|_s^2 \| u - v \|_s + \|v\|_s \| u + v \|_s \| u - v \|_s).
\end{equation}

\end{thm}
{\bf Proof:} 
This result is consequence of the equality
\[
| u|^2 u - |v|^2 v = u^2 ( \bar{u} - \bar{v} ) + \bar{v} (u + v )(u-v )
\]
and the fact that $H_{per}^s$ is a Banach algebra for $s > 1/2$.

\begin{thm}
If $u \in C([0,T]; H_{per}^s)$ with $s > 1/2$ is a solution of (\ref{schrononlinear}), then $u$ satisfies the integral equation
\begin{equation}\label{integraleq}
u(t) = \mathcal{T}(t)u_0 + \int_{0}^{t} \mathcal{T}(t - \xi) F(u(\xi)) d\xi.
\end{equation}
Analogously, if $u \in C([0,T]; H_{per}^s)$, $s > 1/2$ is a solution of (\ref{integraleq}), then $u \in C^1([0,T]; H^s )$ and satisfies (\ref{schrononlinear}) in the following sense:
\begin{equation*}
\lim_{h \to 0^+} \left\|\frac{u(t+h) - u(t)}{h} - \mathcal{A}(u(t)) - F(u(t)) \right\|_{ s } =0
\end{equation*}  
\end{thm}

\textbf{Proof: }  Let $u(t) \in C([0,T]; H^s )$ be a solution of the IVP (\ref{schrononlinear}). Then, for $0 \leq \xi \leq t$ we have
\begin{equation*}
\mathcal{T}(t-\xi) u'(\xi) = \mathcal{T}(t-\xi) \mathcal{A}(u(\xi)) + \mathcal{T}(t-\xi) F(u(\xi)).
\end{equation*}
Note that
\begin{equation*}
\frac{d}{d\xi}\left(\mathcal{T}(t-\xi) u(\xi) \right) = \mathcal{T}(t-\xi) u'(\xi) -\mathcal{T}(t-\xi) \mathcal{A}(u(\xi)).
\end{equation*}
Thus
\begin{equation}\label{derivate}
\frac{d}{d\xi}\left(\mathcal{T}(t-\xi) u(\xi) \right)= \mathcal{T}(t-\xi) F(u(\xi)).
\end{equation}
Integrating at the both sides of the Equation (\ref{derivate}) we obtain that $Y$ satisfies the integral equation:
\begin{equation*}
u(t) = \mathcal{T}(t)u_0 + \int_{0}^{t} \mathcal{T}(t - \xi)  F(u(\xi)) d\xi.
\end{equation*} 
On the other hand, suppose that $u(t) \in C([0,T]; H^s )$ is a solution of the integral equation (\ref{integraleq}). Consider the expression 
\begin{equation}\label{consideration}
\Gamma := \left\| \frac{u(t+h) - u(t)}{h} - \mathcal{A}(u(t)) - F(u(t)) \right\|_{s}.
\end{equation}
By substitution the expressions in (\ref{integraleq}) into (\ref{consideration}) and applying the Minkowsky inequality, we obtain that
\begin{align*}
\Gamma &\leq \left\|\mathcal{T}(t) \left(\frac{\mathcal{T}(h) - I}{h} \right)u_0 - \mathcal{T}(t)\mathcal{A}(Y_0) \right\|_{s } + \\
& \qquad \left\|\left[\frac{\mathcal{T}(h) - I}{h} - \mathcal{A} \right] \int_{0}^{t} \mathcal{T}(t-\xi) F(Y(\xi)) d\xi \right\|_{s} +   \\
&\qquad +\left\|\frac{1}{h} \int_{t}^{t+h} \mathcal{T}(t+h-\xi) F(u(\xi)) d\xi - F(u(t)) \right\|_{s}.
\end{align*}

Note that
\begin{align*}
&\lim_{h \to 0^+} \left\|\mathcal{T}(t) \left(\frac{\mathcal{T}(h) - I}{h} \right)u_0 - \mathcal{T}(t)\mathcal{A}(u_0) \right\|_{s} =0 \\
&\lim_{h \to 0^+} \left\|\left[\frac{\mathcal{T}(h) - I}{h} - \mathcal{A} \right] \int_{0}^{t} \mathcal{T}(t-\xi) F(u(\xi)) d\xi \right\|_{s } = 0.
\end{align*}
The last term can be controlled by using the mean value theorem in the following way:
\begin{align}
&\left\|\frac{1}{h} \int_{t}^{t+h} \mathcal{T}(t+h-\xi) F(u(\xi)) d\xi - F(u(t)) \right\|_{s } \leq \\
& \frac{1}{h} \int_{t}^{t+h} \left\| \mathcal{T}(t+h-\xi) F(u(\xi)) - F(u(t))  \right\|_{s}  d\xi \nonumber  \\
& =\left\|\mathcal{T}(t+h - \tilde{\xi} ) F(u(\tilde{\xi} )) - F(u(t)) \right\|_{s }  \nonumber \\
& \leq \left\| \mathcal{T}(t+h - \tilde{\xi} ) \left[F(u(\tilde{\xi} )) - F(u(t))\right]  \right\|_{s} + \left\|\left[\mathcal{T}(t+h - \tilde{\xi} ) - I \right] F(u(t)) \right\|_{s}, \label{lastterm}
\end{align}
 for some  $t \leq \tilde{\xi} \leq t+h$. Observe that the last term in (\ref{lastterm}) satisfies that
\begin{equation*}
\lim_{h \to 0^+} \left\|\left[\mathcal{T}(t+h - \tilde{\xi} ) - I \right] F(u(t)) \right\|_{s } =0,
\end{equation*}
since $h \to 0^+$ implies $\tilde{\xi} \to 0^+$. Furthermore, for the first term in (\ref{lastterm}) we have
\begin{align*}
& \left\| \mathcal{T}(t+h - \tilde{\xi} )  \left[F(u(\tilde{\xi} )) - F(u(t))\right]  \right\|_{s} \leq \mathcal{K} \left\|  \left[F(u(\tilde{\xi} )) - F(u(t))\right]  \right\|_{s }  \to 0,\\
\end{align*}
as $h \to 0^+$. Hence
\begin{equation*}
\lim_{h \to 0^+} \left\| \frac{u(t+h) - u(t)}{h} - \mathcal{A}(u(t)) - F(u(t)) \right\|_{s} =0.
\end{equation*}

\begin{thm}\label{existencetheorem}
Let $s > 1/2$ and $u_0 \in H_{per}^s$. Then there exists $T>0$ and a unique solution
$u \in C([0, T], H_{per}^s )$ of the integral equation
\[
u(x,t) = \mathcal{T}(t) u_0(x) + \int_0^t \mathcal{T}(t- \xi) F(u(\xi)) d\xi.
\]
\end{thm}

{\bf Proof: }
Let us define the set
\begin{equation*} 
\mathcal{S}_M = \left\{ u  \in C([0,T]; H^s ) : \sup_{t \in [0,T]} \left\|u(t) - \mathcal{T}(t)u_0 \right\|_{s} \leq M   \right\},
\end{equation*}
with the norm $\left\|u \right\|_{s} = \sup_{t \in [0,T]}  ( \left\|u(t) \right\|_{s} )$. Observe that endowed with this norm, $\mathcal{S}_M$ is a  complete set in $C([0,T]; H_{per}^s)$ and $F$ is continuous in $\mathcal{S}_M$.

For $u \in \mathcal{S}_M$ we define the operator
\begin{equation*}
\Psi (u(t)) = \mathcal{T}(t)u_0 + \int_{0}^{t} \mathcal{T}(t-\xi) F(u(\xi)) d\xi.
\end{equation*}
Observe that 
\begin{equation*}
\begin{split}
&\left\| \Psi (u(t+h)) - \Psi (u(t)) \right\|_{s} \\
&\qquad =\Big\| \mathcal{T}(t+h)u_0 + \int_{0}^{t+h} \mathcal{T}(t+h-\xi) F(u(\xi)) d\xi - \mathcal{T}(t)u_0 \\
& \qquad \qquad \qquad \qquad \qquad- \int_{0}^{t} \mathcal{T}(t-\xi) F(u(\xi)) d\xi \Big\|_s \nonumber \\
&\qquad \leq \left\|(\mathcal{T}(t+h) - \mathcal{T}(t)) u_0 \right\|_{s} + \left\| \int_{0}^{t} (\mathcal{T}(t+h-\xi) - \mathcal{T}(t-\xi)) F(u(\xi)) d\xi \right\|_{s} \\
&\qquad \qquad \qquad \qquad \qquad \qquad  \qquad +\left\| \int_{t}^{t+h} \mathcal{T}(t+h-\xi) F(u(\xi)) d\xi\right\|_{s}. 
\end{split}
\end{equation*}
Note that the first term at the right side above satisfies
\begin{equation*}
\lim_{h \to 0^+} \left\|(\mathcal{T}(t+h) - \mathcal{T}(t)) u_0 \right\|_{s } =0,
\end{equation*}
since $(\mathcal{T}(t))_{t \geq 0}$ is a $C_0$-semigroup in $H_{per}^s$. The second and third terms at the right side are controlled using the Lebesgue's dominated convergence theorem. Note that
\begin{equation*}
\left\| (\mathcal{T}(t+h-\xi) - \mathcal{T}(t-\xi)) F(u(\xi)) \right\|_{s} \leq \mathcal{K} \left\| F(Y(\xi)) \right\|_{s } \leq \mathcal{K} \left\| F(u(\xi)) \right\|_{s },
\end{equation*}
and
\begin{equation*}
\lim_{h \to 0^+} \left\| (\mathcal{T}(t+h-\xi) - \mathcal{T}(t-\xi)) F(u(\xi)) \right\|_{s} = 0.
\end{equation*}
Thus,
\begin{align*}
\lim_{h \to 0^+} & \left\| \int_{0}^{t} (\mathcal{T}(t+h-\xi) - \mathcal{T}(t-\xi)) F(u(\xi)) d\xi \right\|_{s }   \\
 & \leq \lim_{h \to 0^+} \int_{0}^{t}  \Big\|  (\mathcal{T}(t+h-\xi) - \mathcal{T}(t-\xi)) F(u(\xi)) \Big\|_{s } d\xi =0.
\end{align*}
Finally,
\begin{align*}
\left\|\mathcal{T}(t+h-\xi) F(u(\xi)) \right\|_{s } &\leq \mathcal{K} \left\|u(\xi) \right\|_{s} \\
&= \mathcal{K} (\left\| u \right\|_{s}  )(\left\| u \right\|_{s} ) \\
&= \mathcal{K} (\left\|u \right\|_{s} )^2 \\
&= \mathcal{K} \left\|u \right\|_{s}^2 \\
&\leq (M+ \mathcal{K} \left\| u_0 \right\|_{s} )^2,
\end{align*}
which is a consequence of the estimate
\begin{align*}
\left\|u(t) \right\|_{s} &= \left\| u(t) - \mathcal{T}(t) u_0 + \mathcal{T}(t) u_0 \right\|_{s} \\
&\leq \left\| u(t) - \mathcal{T}(t)u_0 \right\|_{s}  + \left\| \mathcal{T}(t) u_0 \right\|_{s} \\
&\leq (M +  \mathcal{K}\left\|u_0 \right\|_{s}).
\end{align*}
Thus,
\begin{equation*}
\lim_{h \to 0^+} \left\|\mathcal{T}(t+h-\xi) F(u(\xi)) \right\|_{s } =0.
\end{equation*}
It follows that
\begin{align*}
\lim_{h \to 0^+}  \left\| \int_{t}^{t+h} \mathcal{T}(t+h-\xi) F(u(\xi)) d\xi \right\|_{s } & \leq \lim_{h \to 0^+}  \int_{t}^{t+h} \left\|  \mathcal{T}(t+h-\xi) F(u(\xi)) \right\|_{s } d\xi \\
& =0.
\end{align*}
Hence, if $u \in \mathcal{S}_M$, then $\Psi (u) \in C([0,T]; H^s)$. 

On the other hand, 
\begin{align*}
\left\|\Psi (u(t)) - \mathcal{T}(t)u_0 \right\|_{s } &= \left\|\int_0^t \mathcal{T}(t-\xi) F(u(\xi)) d\xi \right\|_{s } \\
&\leq \int_{0}^{t} L(\left\|u \right\|_{s},0) \left\|u(\xi) \right\|_{s} d\xi. \\
&\leq T \mathcal{K} L(M +  \mathcal{K}\left\| u_0 \right\|_{s}, 0) (M +  \mathcal{K}\left\|u_0 \right\|_{s} ).
\end{align*}
Here $L(.,.)$ is a polynomial.  Thus choosing
\begin{equation*}
T_1 = \frac{M}{\mathcal{K} L(M +  \mathcal{K}\left\|u_0 \right\|_{s},0,M +  \mathcal{K}\left\|u_0 \right\|_{s}) (M +  \mathcal{K}\left\| u_0 \right\|_{s})},
\end{equation*}
we obtain that $\Psi(u(t)) \in \mathcal{S}_M$ if $u(t) \in \mathcal{S}_M$.

\

To see that there exists $T_2$ such that $\Psi$ is a contraction in $\mathcal{S}_M$ for $T < T_2$, let us observe that
\begin{align*}
\left\| \Psi(u(t)) - \Psi( \overline{u}(t)) \right\|_{s } &= \left\| \int_{0}^{t} \mathcal{T}(t-\xi)\left[F(u(\xi)) - F(\overline{u}(\xi))\right] d\xi \right\|_{s }  \\
&\leq \mathcal{K} T L(M +  \mathcal{K}\left\|u_0 \right\|_{s},0 ) \sup_{\xi \in [0,T]}  \left\| u(\xi) - \overline{u}(\xi) \right\|_{s}. 
\end{align*}
Thus choosing
\begin{equation*}
T_2 = \frac{1}{\mathcal{K} L(M +  \mathcal{K}\left\|u_0 \right\|_{s},0 )},
\end{equation*}
we obtain that $\Psi$ is a contraction on $\mathcal{S}_M$. Thus with $T^* < \min \{T_1,T_2 \}$ and by applying Banach's fixed-point theorem on $\mathcal{S}_M$, we get the desired result.

\begin{thm}
The solutions obtained in Theorem \ref{existencetheorem} is unique and depends continuously on the initial condition $u_0$.
\end{thm}

\textbf{Proof:} Let $u = $ and $\overline{u}$ be elements in $C([0,T]; H_{per}^s)$ solutions of the integral equation (\ref{integraleq}) with initial data $u_0$ and $\overline{u}_0$, respectively. Then
\begin{align*}
\left\| u(t) - \overline{u}(t) \right\|_{s } &= \left\| \mathcal{T}(t) (u_0 - \overline{u}_0) + \int_{0}^{t} \mathcal{T}(t-\xi) \left[F(u(\xi)) - F(\overline{u}(\xi)) \right] d\xi \right\|_{s } \\
&\leq \mathcal{K} \left\| u_0 - \overline{u}_0 \right\|_{s } + \int_{0}^{t} \mathcal{K} L(\left\|u \right\|_{s},\left\|\overline{u} \right\|_{s} ) \left\| u(\xi) - \overline{u}(\xi) \right\|_{s} d\xi \\
&\leq \mathcal{K} \left(\left\|u_0 - \overline{u}_0 \right\|_{s } + L(\left\|u \right\|_{s},\left\|\overline{u} \right\|_{s} )\int_{0}^{t}   \left\|u(\xi) - \overline{u}(\xi) \right\|_{s} d\xi \right).
\end{align*}
Let $\overline{L} := \sup_{t \in [0,T]} L(\left\|u(t) \right\|_{s}, \left\|\overline{u}(t) \right\|_{s} )$. The Gronwall's inequality implies that
\begin{equation*}
\left\| u(t) - \overline{u}(t) \right\|_{s } \leq \mathcal{K} \left\| u_0 - \overline{u}_0 \right\|_{s}  e^{\overline{L} t}) \leq \mathcal{K} \left\| u_0 - \overline{u}_0 \right\|_{s }  e^{\overline{L} T},
\end{equation*}
for all $t \in [0,T]$. Hence uniqueness and continuous dependence of the solutions on initial conditions follows.


\section{Study of the limits $\eta \to 0^+$, $\beta \to 0^+$}

In this section, we analyze the limits $\eta \to 0^+$ and $\gamma \to 0^+$ in the equation
\begin{equation}\label{full_model}
\partial_t u - i \partial_{x}^2 u + \alpha \partial_x ( |u|^2 u ) = \eta \partial_x^2 u + \beta \partial_x^3 u + \gamma \partial_x u.
\end{equation}

The following Kato's commutator lemma will be important to estimate the derivative nonlinear term in the
previous equation.
\begin{lemma} \label{Kato_estimate} \cite{Kato_commutator}
Let $s >3/2$, and $t \geq 1$. If f, u are real valued, then
\begin{equation}\label{Kato_estimate}
| \langle f Du, u \rangle_t | \leq C (\| D f \|_{s-1} \| u \|_t^2 + \| D f \|_{t-1} \| u \|_s \| u \|_t),
\end{equation}
where $D$ is the first-order derivative.
\end{lemma}

\begin{thm}\label{uniform_bound}
Let $\eta >0$, $s> 3/2$, $u_0 \in H_{per}^s$, and $u = u_\eta$ be the solution of problem
\eqref{schrononlinear} with $u_\eta(0) = u_0$. Then there exists a constant $C>0$ independent of
$\eta$, such as for any $0 < T_\eta < \frac{2}{C \|u_0\|_s}$, the solution $u$ can be extended to the interval $[0, T_s]$, and the estimate
\[
\| u_\eta \|_s \leq \frac{2 \| u_0\|_s }{2 - C t \| u_0 \|_s}
\]
holds.
\end{thm}

{\bf Proof:} Let $u = u_\eta$ be a solution of equation \eqref{schrononlinear} with $u_\eta(0) = u_0$.
Multiplying  this equation by $u$, we obtain
\[
\langle \partial_t u, u \rangle_s - i \langle \partial_x^2 u , u \rangle_s + \alpha \langle \partial_x(|u|^2 u), u \rangle_s
= \langle \eta \partial_x^2 u, u \rangle_s + \beta \langle \partial_x^3 u, u \rangle_s + \gamma \langle \partial_x u, u \rangle_s.
\]
But $\langle \partial_x^3 u, u \rangle_s = \langle \partial_x u, u \rangle_s = 0$, and thus,
\begin{equation}\label{norm_u}
\langle \partial_t u, u \rangle_s + i \| \partial_x u \|_s^2 + \eta \| \partial_x^2 u \|_s^2 = - \alpha \langle \partial_x (|u|^2 u ), u \rangle_s.
\end{equation}
Now since
\[
\frac{d}{dt} \|u\|_s^2 = 2 \text{Re} \langle \partial_t u, u \rangle_s,
\]
we obtain after taking the real part of both sides of equation \eqref{norm_u},
\begin{align}
\frac{1}{2} \frac{d}{dt} \| u \|_s^2 + \eta \| \partial_x^2 u \|_s^2 &= - \alpha \text{Re} \langle \partial_x (|u|^2 u ), u \rangle_s\\
& \leq |\alpha \langle 2|u|^2 \partial_x u , u \rangle_s | + | \alpha \langle u^2 \partial_x \bar{u} , u \rangle_s |
\end{align}
Now using the Kato's commutator estimate \eqref{Kato_estimate} for the nonlinear terms on the right hand side of the previous equation, and the fact that
\[
\| u \|_{s}^2 = \| \text{Re} (u) \|_{s}^2 + \| \text{Im}(u) \|_s^2,
\]
we obtain
\[
\frac{1}{2} \frac{d}{dt} \| u \|_s^2 + \eta \| \partial_x^2 u \|_s^2  \leq C  \| u \|_s^4,
\]
where $C>0$ is a constant independent of $u$. As a consequence, since $\eta >0$,
\begin{equation}
\frac{d}{dt} \| u \|_s^2   \leq C  \| u \|_s^4,
\end{equation}
or equivalently,
\begin{equation}
\frac{d}{dt} z(t) \leq C z(t)^2,
\end{equation}
where $z(t) = \|u\|_s^2$.

Integrating in $[0,t]$,
\[
\frac{1}{z(t)} \geq \frac{1}{z(0)} - \frac{C t}{2}. 
\]
Since $z(0) = \| u(0) \|_s^2 = \| u_0 \|_s^2$, we obtain
\begin{equation}\label{solution_estimate}
\| u(t) \|_s^2 = z(t) = \frac{2 \| u_0 \|_s^2 }{ 2 - C \|u_0\|_s^2 t},
\end{equation}
whenever $t < T_s < \frac{2}{C \|u_0\|_s^2}$.

If $T_\eta < T_s$, we can solve the problem \eqref{schrononlinear} for $t \geq T_\eta$, with initial value $u_\eta(T_\eta) \in H_{per}^s$, obtaining a solution defined in the interval $[0, T_\eta + T_\eta']$ (by using Theorem \ref{existencetheorem} ) and such that the inequality \eqref{solution_estimate} is satisfied. Given that $T_\eta' >0$ can be determined depending only on $\eta$ and $u_\eta(T_\eta)$, we have that in a finite number of such extensions the solution can be extended to a solution defined in $[0, T_s]$ where the estimate \eqref{solution_estimate} is satisfied.

\begin{thm}
Let $s > 3/2$, $u_0 \in H_{per}^s$. Then the solution $u_\eta$ of problem \eqref{full_model} subject to $u_\eta(0) = u_0$ converges weakly in $H_{per}^s$ as $\eta \to 0^+$ to an integral solution $u \in C([0, T_s], L_{per}^2 )$ of the problem

\begin{align} \label{schro2}
i \partial_t u + \partial_x^2 u + i \alpha \partial_x ( |u|^2 u ) =  
i \beta \partial_x^3 u + i \gamma \partial_x u,
\end{align}
subject to the initial condition 
\begin{equation}\label{initial_cond2}
u(x,0)=u_0(x).
\end{equation}

\end{thm}

{\bf Proof: } 

Let $u_\epsilon, u_\delta \in C([0, T_s]; H_{per}^s )$ be solutions of problem \eqref{full_model} with $\eta = \epsilon>0, \eta=\delta >0$. Here $0< T_s < \frac{2}{C \|u_0\|_s^2}$ as in Theorem \ref{uniform_bound}. Let us define $w = u_\epsilon-u_\delta$. Therefore $w$ satisfies the equation
\[
\partial_t w - i \partial_x^2 w + \alpha \partial_x( |u_\epsilon|^2 u_\epsilon - |u_\delta|^2 u_\delta ) -\beta \partial_x^3 w - \gamma \partial_x w = \partial_x^2 ( \epsilon u_\epsilon - \delta u_\delta ).
\]
By taking the $L^2$ inner product of the previous equation with $w$, and using that
$\langle \partial_x^3 w, w \rangle = \langle \partial_x w, w \rangle = 0$, we have
\[
\langle \partial_t w, w \rangle - i \langle \partial_x^2 w , w \rangle = - \alpha \langle \partial_x ( |u_\epsilon|^2 u_\epsilon - |u_\delta |^2 u_\delta) , w \rangle + \epsilon \langle \partial_x^2 w , w \rangle + (\epsilon- \delta ) \langle \partial_x^2 u_\delta, w \rangle.
\]
Taking real part of the previous equation, we have
\begin{equation}\label{estimate_w}
\frac12 \frac{d}{dt} \|w\|^2 + \epsilon \| \partial_x w \|^2 = -\alpha \text{Re} \langle \partial_x ( |u_\epsilon|^2 u_\epsilon - |u_\delta |^2 u_\delta), w \rangle +(\epsilon-\delta ) \text{Re} \langle \partial_x^2 u_\delta, w \rangle.
\end{equation}

\noindent Now by virtue of the identity,
\[
|f|^2 f - |g|^2 g = |f|^2 (f - g) + fg (\overline{f} - \overline{g}) + |g|^2 (f - g),
\]
we obtain using integration by parts,
\[
-\alpha \text{Re} \langle \partial_x( |u_\epsilon|^2 u_\epsilon - | u_\delta |^2 u_\delta ), w \rangle = \alpha \int_0^L \text{Re}( ( | u_\epsilon |^2 w + u_\epsilon u_\delta \overline{w} + |u_\delta|^2 w ) \partial_x \overline{w} )dx
\]
Since
\[
\text{Re} ( w \partial_x \overline{w} ) = \frac12 \partial_x( w \overline{w} ), ~~~~ \ \overline{w} \partial_x \overline{w}  = \frac12 \partial_x(  \overline{w}^2 ),
\]
we obtain using the continuous embedding $H_{per}^{s-1} \subset L_{per}^\infty$, for $s >3/2$,
\begin{align*}
&-\alpha \text{Re} \langle \partial_x( |u_\epsilon|^2 u_\epsilon - | u_\delta |^2 u_\delta ), w \rangle = \frac{\alpha}{2}  \int_0^L  \Big( |u_\epsilon|^2  \partial_x( w \overline{w} ) + u_\epsilon u_\delta \partial_x( \overline{w}^2) + |u_\delta |^2 \partial_x ( w \overline{w} ) \Big)  dx\\
&=-\frac{\alpha}{2} \int_0^L \Big(  \partial_x( |u_\epsilon|^2 ) w \overline{w} + \partial_x(u_\epsilon u_\delta) \overline{w}^2 + \partial_x(| u_\delta|^2) w \overline{w} \Big) dx\\
&=-\frac{\alpha}{2} \int_0^L \Big(  \partial_x(|u_\epsilon|^2) |w|^2 + \partial_x(u_\epsilon u_\delta) \overline{w}^2 + \partial_x( | u_\delta|^2 ) |w|^2 \Big) dx\\
& \leq C \Big( \| \partial_x(|u_\epsilon|^2) \|_\infty \| w \|^2 + \| \partial_x( u_\epsilon u_\delta) \|_\infty \|w\|^2 + \|\partial_x(|u_\delta|^2 )  \|_\infty  \|w\|^2 \Big) \\
&\leq C \Big( \| \partial_x(|u_\epsilon|^2) \|_{s-1} + \| \partial_x( u_\epsilon u_\delta) \|_{s-1}  + \|\partial_x(|u_\delta|^2 )  \|_{s-1}   \Big) \|w \|^2 \\
&\leq C \Big( \| |u_\epsilon|^2 \|_s  + \| u_\epsilon u_\delta \|_s  + \| |u_\delta|^2 \|_s   \Big) \| w \|^2.
\end{align*}
From Theorem \ref{uniform_bound}, the sequence $(u_\epsilon)_{\epsilon>0}$ is uniformly bounded in $[0, T_s]$.
Moreover, from the estimate \eqref{estimate_w}, we have
\[
\frac12 \frac{d}{dt} \| w \|^2 \leq C \| w \|^2 + C | \epsilon - \delta |, ~~ 0 \leq t \leq T_s,
\]
where $C$ is a constant independent of $\epsilon$ and $\delta$. By applying Gronwall's inequality, it follows that
\[
\| w(t) \| ^2 \leq T_s C e^{C T_s} | \epsilon - \delta |, ~~~ 0 \leq t \leq T_s.
\]
Therefore, $(u_\epsilon)_{\epsilon >0}$ is a Cauchy sequence in the Banach space
normed space
\[
C([0,T_s]; L_{per}^2),
\]
and hence converges strongly as $\epsilon \to 0^+$ to some function $u \in C([0,T_s]; L_{per}^2)$, uniformly for $t \in [0, T_s]$. 

Since $(u_\epsilon)$ is uniformly bounded in $L^\infty([0,T_s]; H_{per}^s)$, it follows that the sequence $(u_\epsilon)$ converges weakly in $H_{per}^s$ to $u$, uniformly for $t \in [0,T_s]$. To see this, let $\psi \in H_{per}^s$ and $\varphi \in \mathcal{C}_{per}^\infty(0,L):=S_{per}$ (space of infinitely differentiable $L$-periodic functions).  Then, using the Cauchy-Schwarz inequality, we have
\begin{align*}
|\langle u_\epsilon-u_\delta, \psi \rangle_s | &= |\langle u_\epsilon - u_\delta , \psi + \varphi -\varphi \rangle_s | \\
&\leq | \langle u_\epsilon - u_\delta , \psi - \varphi \rangle_s | + | \langle u_\epsilon - u_\delta ,\varphi \rangle_s | \\
&\| u_\epsilon - u_\delta \|_s \| \psi - \varphi \|_s + \| u_\epsilon - u_\delta \| \| \varphi \|_{2s}.
\end{align*}
Since $S_{per}$ is dense in $H_{per}^s$, and $u_\epsilon$ converges strongly in $L_{per}^2$, it follows that $(u_\epsilon)$ is a weakly Cauchy sequence in $H_{per}^s$. Hence, we conclude that $u_\epsilon(t) \rightharpoonup u(t)$ in $H_{per}^s$, uniformly  for $t \in [0, T_s]$.

Now consider the nonlinear operator
\[
F_\epsilon (u) := i \partial_x^2 u - \alpha \partial_x (|u|^2 u ) + \epsilon \partial_x^2 u + \beta \partial_x^3 u + \gamma \partial_x u.
\]
The function $u_\epsilon$ satisfies the differential equation
\begin{equation}\label{EDO_eq}
\frac{d}{dt} u_\epsilon(t) = F_\epsilon( u_\epsilon(t) ).
\end{equation}
It can be shown that 
\[
F_\epsilon (u_\epsilon) \rightharpoonup  F_0(u),
\]
as $\epsilon \to 0^+$ in $H_{per}^{s-3}$. 

Integrating equation \eqref{EDO_eq} over the interval $[\tau, t]$, we obtain
\[
u_\epsilon(t) = u_\epsilon(\tau) + \int_\tau^t F_\epsilon (u_\epsilon(t') ) dt'.
\]
Passing to the limit $\epsilon \to 0^+$, and using the previous results, we obtain
\begin{equation}\label{integral_sol}
u(t) = u(\tau) + \int_\tau^t F_0 (u(t') ) dt'.
\end{equation}

We now show that $\lim_{t \to 0^+} u(t) = u_0$ in $H_{per}^s$.  Let us see the continuity from the right of $u$ at $t=0$. Observe that selecting $\psi \in H_{per}^s$ and $\epsilon>0$ small enough, we have that
\begin{align*}
| \langle u(t) , \psi \rangle_s - \langle u_0, \psi \rangle_s | & \leq | \langle u(t) , \psi \rangle_s - \langle u_\epsilon(t), \psi \rangle_s | + | \langle u_\epsilon(t) , \psi \rangle_s - \langle u_0, \psi \rangle_s | \\
&\leq | \langle u_\epsilon(t) - u(t), \psi \rangle_s | + | \langle u_\epsilon(t) - u_0, \psi \rangle_s |.
\end{align*}
From the weak convergence and  $u_\epsilon \in C([0, T_s]; H_{per}^s)$, follows that $u_\epsilon(t)$ converges weakly to $u_0$ in $H^s$, when $t \to 0^+$, with $t \in [0, T_s]$. As a consequence, 
\[
\| u_0 \|_s \leq \lim_{t \to 0^+} \text{inf} \| u(t) \|_s.
\]
Now, due to estimate \eqref{solution_estimate}, we have that
\[
\lim_{t \to 0^+} \text{sup} \|u(t)\|_s^2 \leq \lim_{t \to 0^+} \text{sup} \frac{2 \|u_0\|_s^2}{2 - C \|u_0\|_s^2 t} = \| u_0 \|_s^2.
\]
Therefore,
\[
\|u_0\|_s \leq \lim_{t \to 0^+} \text{inf} \| u(t)\|_s \leq \lim_{t \to 0^+} \text{sup} \| u(t)\|_s \leq \|u_0\|_s.
\]
In this manner, we conclude that
\begin{equation}\label{convergence_norm}
\lim_{t \to 0^+} \|u(t)\|_s = \|u_0\|_s.
\end{equation}
But
\begin{align*}
\| u(t) - u_0\|_s &= \langle u - u_0, u- u_0 \rangle_s = \langle u, u \rangle_s - \langle u_0, u\rangle_s  - \langle u, u_0 \rangle_s + \langle u_0, u_0 \rangle_s\\
&\leq | \langle u, u \rangle_s - \langle u_0, u_0 \rangle_s | + | \langle u_0, u_0 \rangle_s  - \langle u, u_0 \rangle_s | + | \langle u_0,u_0 \rangle_s - \langle u_0 , u \rangle_s | \\
&\leq | \|u\|_s^2 - \|u_0\|_s^2 | + | \langle u_0 - u, u_0 \rangle_s | + | \langle u_0, u_0 - u \rangle_s |.
\end{align*}
Thus from \eqref{convergence_norm} and the fact that $u(t) \rightharpoonup  u_0$, $t \to 0^+$, follows that
$\lim_{t \to 0^+} u(t) = u_0$ in $H_{per}^s$.

Taking limit $\tau \to 0^+$ in equation \eqref{integral_sol}, we obtain
\[
u(t) = u_0 + \int_0^t F_0(u(t')) dt',
\]
which shows that $u$ is an integral solution of equation \eqref{full_model}.


\bigskip

\begin{thm}
Let $s > 3/2$, $u_0 \in H_{per}^s$. Then the solution $u_\beta$ of problem \eqref{full_model} subject to $u_\beta(0) = u_0$ converges weakly in $H_{per}^s$ as $\beta \to 0^+$ to an integral solution $u \in C([0, T_s], L_{per}^2)$ of the problem
\begin{align} \label{schro3}
i \partial_t u + \partial_x^2 u + i \alpha \partial_x ( |u|^2 u ) =  
i  \eta \partial_x^2 u + i \gamma \partial_x u, 
\end{align}
subject to the initial condition 
\begin{equation}\label{initial_cond3}
u(x,0)=u_0(x).
\end{equation}

\end{thm}

{\bf Proof: } The proof is analogous to the previous theorem.

\section{The numerical scheme}
 
To approximate the solutions of the periodic initial value problem \eqref{schro1}-\eqref{periodic_cond}, the interval $[0,L]$ is discretized by $N$ equally spaced points with distance $\Delta x= L/N$. Let us consider the Fourier basis 
$$
\phi_k(x) = e^{i \omega_k x}, ~~ \omega_k = \frac{2\pi k}{L}, ~~~ k= -N/2,...,N/2.
$$
It is well known that the set $\{ L^{-1/2} \phi_k \}_k$ is an orthonormal basis of the space $L^2(0,L)$, with respect to the inner product
$$
\langle f,g \rangle =  \int_0^L f(x) \overline{ g(x) } d x.
$$
For $N \in 2 \mathbb{N}$, let us consider the finite dimensional space
$$S_N:= \text{span} \{ L^{-1/2} \phi_k, k=-N/2,...N/2 \}.$$ 
We denote by $P_N$ the orthogonal projection operator
 $P_N : L_{per}^2(0,L) \to S_N$ defined by
$$
P_N g = \sum_{k=-N/2}^{N/2} \hat{g}_k \phi_k,
$$
where
$$
\hat{g}_k = \frac{1}{L} \int_0^L g(x) \overline{\phi}_k(x) dx,~~~ k \in \mathbb{Z}.
$$
Since the linear operator $P_N$ is a orthogonal projection, we have that
$$
\langle P_N g -g, \phi \rangle = 0,
$$
for all $\phi \in S_N$.

Furthermore, given real numbres $r,s$, with $0 \leq s \leq r$, there exists a constant $C>0$ such
that, for any $u \in H_{per}^r(0,L)$, 
\begin{equation}\label{projection_prop}
\| P_N u - u \|_s \leq C N^{s-r} \| u \|_r.
\end{equation}

Thus, the weak formulation of equation \eqref{schro1} consists in finding $u \in C^1(0,T; H^2(0,L))$ such that
\begin{equation} \label{weak}
\begin{aligned}
&  i \partial_t \langle u, \phi \rangle -\langle \partial_x u, \partial_x \phi \rangle -i \alpha \langle |u|^2 u, \partial_x \phi \rangle = -i \langle \eta \partial_x u + \beta \partial_x^2 u + \gamma u, \partial_x \phi \rangle,  \\
& u(x,0) = u_0(x), 
\end{aligned}
\end{equation}
for all $\phi \in H^1(0,L)$, $t \in [0,T]$.

\subsection{Semidiscrete scheme}

The semidiscrete Galerkin spectral method to approximate problem (\ref{weak}) based on the Fourier basis $\phi_k$ consists in finding $u_N \in C^1(0,T_N; S_N)$ such that
\begin{align}
&i \partial_t \langle u_N, \phi \rangle  - \langle \partial_x u_N, \partial_x \phi \rangle -i \alpha \langle |u_N|^2 u_N, \partial_x \phi \rangle  = \notag \\ 
& \, \, \, \, \, \, \, \, \, \, \, \, \, \,  -i \langle \eta \partial_x u_N + \beta \partial_x^2 u_N + \gamma u_N, \partial_x \phi \rangle, \label{semidis} \\
&u_N(0) = P_N(u_0), \notag
\end{align}
for all $\phi \in S_N$, $t \in [0,T_N]$. Observe that since $u_N \in S_N$, we can write
\begin{equation}
u_N(x,t) = \sum_{k=-N/2}^{N/2} \hat{u}_k(t) \phi_k(x).  \label{approximation}
\end{equation}
Substituting $\phi = L^{-1/2} \phi_k$ in equation \eqref{semidis} and using the expansion \eqref{approximation}, we obtain
\begin{equation}\label{semidis2}
i \hat{u}'_k(t) + ( -w_k^2 + i \eta w_k^2 - \beta w_k^3 + \gamma w_k) \hat{u}_k =\alpha w_k P_k [ |u_N|^2 u_N ],
\end{equation}
where
\[
P_k[g] =  \frac{1}{L} \int_0^L g(x) \overline{\phi}_k(x) dx.
\]

The following theorem establishes the spectral accuracy of the proposed numerical solver \eqref{semidis}:

\begin{thm}
Let $s > 3/2$ be an integer and $u \in C([0,T], H_{per}^r)$ a classical solution of equation \eqref{schrononlinear}, corresponding to initial data $u_0 \in H_{per}^r$ for some $r > s$. Then the semidiscrete problem
\eqref{semidis} has a unique solution $u_N \in C([0,T], S_N)$, which satisfies for $N$ large enough and some constant $C>0$ independent of $N, t$, such that
\[
\|u(t) - u_N(t) \|_s \leq C N^{s - r},
\]
for any $0 \leq t \leq T$.

\end{thm}

{\bf Proof: }  From the classical theory of ordinary differential equations, one can show that the initial value problem \eqref{semidis} has a unique solution $u_N \in C([0,T_N], S_N)$.  

Notice that from equation \eqref{schro1}  and \eqref{semidis}, we get for any $\phi \in S_N$ that
\begin{align*}
&i \langle P_N \partial_t u - \partial_t u_N, \phi \rangle_s + \langle P_N \partial_x^2 u - \partial_x^2 u_N, \phi \rangle_s + i\alpha \langle P_N(\partial_x (|u|^2 u) ) - \partial_x (|u_N|^2 u_N), \phi \rangle_s\\
&-i \eta \langle P_N (\partial_x^2 u ) - \partial_x^2 u_N, \phi \rangle_s - i \beta \langle P_N (\partial_x^3 u ) - \partial_x^3 u_N , \phi \rangle_s - \\ 
&i \gamma \langle P_N \partial_x u - \partial_x u_N, \phi \rangle_s = 0.
\end{align*}
Letting 
\[
\theta := P_N u - u_N, ~~~~ \phi = \theta,
\]
we have that
\begin{align*}
&i \langle \partial_t \theta, \theta \rangle_s + \langle \partial_x^2 \theta, \theta \rangle_s + i\alpha \langle
P_N ( \partial_x (|u|^2 u) - \partial_x (|u_N|^2 u_N ), \theta \rangle_s - i \eta \langle \partial_x^2 \theta, \theta \rangle_s - \\
& i \beta \langle \partial_x^3 \theta, \theta \rangle_s - i \gamma \langle \partial_x \theta, \theta \rangle_s = 0.
\end{align*}
But since $\langle \partial_x^3 \theta, \theta \rangle_s = \langle \partial_x \theta , \theta \rangle_s=0$,
and multiplying by $-i$, we have 
\[
\langle \partial_t \theta, \theta \rangle_s + i \| \partial_x^2 \theta \|_s^2 + \alpha \langle
P_N ( \partial_x (|u|^2 u) ) - \partial_x (|u_N|^2 u_N ), \theta \rangle_s +  \eta \| \partial_x^2 \theta \|_s^2 = 0.
\]
Now, since
\[
\text{Re} \langle \partial_t \theta , \theta \rangle_s = \frac12 \frac{d}{dt} \|\theta\|_s^2,
\]
and taking real part on both sides of the previous equality, we arrive at
\begin{align*}
&\frac12 \frac{d}{dt} \| \theta \|_s^2 + \eta \| \partial_x \theta \|_s^2 = - \alpha  \text{Re} [
\langle P_N ( \partial_x (|u|^2 u) - \partial_x (|u_N|^2 u_N ), \theta \rangle_s ] \\
&= -\alpha \text{Re}[ \langle \partial_x (|u|^2 u) - \partial_x (|u_N|^2 u_N ), \theta \rangle_s ]\\
&\leq |\alpha| | \langle |u|^2 u - |u_N|^2 u_N , \partial_x \theta \rangle_s |.
\end{align*}
Now, since 
\[
\| |u|^2 u - |u_N|^2 u_N \|_s \leq C ( \| u \|_s^2 + \|u\|_s \|u_N\|_s + \| u_N \|_s^2 ) \| u - u_N \|_s,
\]
we obtain that
\[
\frac12 \frac{d}{dt} \| \theta \|_s^2 + \eta \| \partial_x^2 \theta \|_s^2 \leq C ( \| u\|_s^2 + \|u\|_s \|u_N\|_s + \|u_N\|_s^2 ) \| u - u_N \|_s \| \partial_x \theta\|_s.
\]
Next, using the H\"older inequality, we have for $\epsilon >0$,
\begin{align*}
&\frac12 \frac{d}{dt} \| \theta \|_s^2 + \eta \| \partial_x^2 \theta \|_s^2 \leq \frac{C}{2 \epsilon} ( \| u\|_s^2 + \|u\|_s \|u_N\|_s + \|u_N\|_s^2 )^2 \| u - u_N \|_s^2 + \epsilon  \| \partial_x \theta \|_s^2.
\end{align*}
Taking into account that $\| \theta \|_s = \|u - u_N\|_s$, we get
\begin{align*}
&\frac12 \frac{d}{dt} \| \theta \|_s^2 + (\eta - \epsilon) \| \partial_x^2 \theta \|_s^2 \leq \frac{C}{2 \epsilon} ( \| u\|_s^2 + \|u\|_s \|u_N\|_s + \|u_N\|_s^2 )^2 \| u - u_N \|_s \|\theta \|_s.
\end{align*}
By selecting $0< \epsilon < \eta$, we obtain
\begin{align*}
&\frac12 \frac{d}{dt} \| \theta \|_s^2 = \| \theta\|_s \frac{d}{dt} \| \theta \|_s  \leq \frac{C}{2 \epsilon} ( \| u\|_s^2 + \|u\|_s \|u_N\|_s + \|u_N\|_s^2 )^2 \| u - u_N \|_s \|\theta \|_s.
\end{align*}
Therefore,
\begin{align}\label{semiestimate}
\frac{d}{dt}  \| \theta \|_s  \leq \frac{C}{2 \epsilon} ( \| u\|_s^2 + \|u\|_s \|u_N\|_s + \|u_N\|_s^2 )^2 \| u - u_N \|_s.
\end{align}
Suppose that for some constant $B>0$,
\[
\| u(t) \|_s \leq B,
\]
for any $0 \leq t \leq T$, and let $0 < T_N < T$ be the maximal time for which
\begin{equation}\label{estimate1}
\| u_N(t) \|_s \leq 2B,
\end{equation}
for $0 \leq t \leq T_N$. We get from inequality \eqref{semiestimate} for some constant $C>0$ independent on $N$ that
\[
\frac{d}{dt} \| \theta \|_s \leq C \|u - u_N \|_s.
\]
Therefore for $r > s$,
\[
\frac{d}{dt} \| \theta \|_s \leq C \|u - u_N \|_s \leq C( \|u - P_N u\|_s + \| P_N u - u_N \|_s \leq C N^{s-r} + \| \theta \|_s.
\]
Let us note that from $\theta(0) = 0$, the last inequality and Gronwall's lemma, we obtain for $0 \leq t \leq T_N$ that
\begin{equation}\label{estimate2}
\| \theta \|_s \leq C N^{s-r},
\end{equation}
where $C>0$ is a constant independent on. $N$.

\noindent On the other hand, since $u_N = u - (u - P_N u + \theta )$, one gets that
\[
\| u_N(t) \|_s \leq \| u(t) \|_s + \| u(t) - P_N u(t) \|_s + \| \theta(t) \|_s \leq B + C N^{s-r}, ~~ 0 \leq t \leq T_N.
\]
Since $r > s$, we conclude that for $N$ large enough
\begin{equation}\label{estimateu_N}
\| u_N(t) \|_s < 2B, ~~~~ 0 \leq t \leq T_N.
\end{equation}
For this reason, the numerical solution $u_N$ can be extended satisfying \eqref{estimateu_N} for
$t \in [0, T_N + \delta]$, for some $\delta >0$. This fact contradicts the maximality of $T_N$. Thus the inequalities \eqref{estimate1}, \eqref{estimate2} are satisfied until $t = T$ and
\[
\| u(t) - u_N(t)\|_s \leq \| u(t) - P_N u(t) \|_s + \| \theta \|_s \leq C N^{s-r},
\]
for any $0 \leq t \leq T$, which finishes the proof of the theorem.

\subsection{Fully discrete scheme}

The semidiscrete formulation given in \eqref{semidis2} is approximated by using the
following second-order finite difference scheme in time:
\begin{align}
&i \Big( \frac{\hat{u}_k^{(n+1)} - \hat{u}_k^{(n)} }{\Delta t}\Big) + ( -w_k^2 + i \eta w_k^2 -\beta w_k^3 + \gamma w_k ) \Big(  \frac{\hat{u}_k^{(n+1)} +  \hat{u}_k^{(n)} }{2} \Big) \label{fullscheme} \\
&~~~~~ = \frac{3}{2} \alpha w_k P_k[ |u_N|^2 u_N ]^{(n)} - \frac{1}{2} \alpha w_k P_k[ |u_N|^2 u_N ]^{(n-1)}, \notag
\end{align}
where $\Delta t$ denotes the time step, and $\hat{u}_k^{(n)}$ represents the numerical
aproximation of the Fourier coefficient $\hat{u}_k(t)$ at time $t = n \Delta t$. The notation $P_j[g]^{(n)}$ indicates the value of $P_k[g]$ when $g$ is evaluated at time $t = n\Delta t$.

\begin{figure}[h]
\includegraphics[width=15cm, height=10cm]{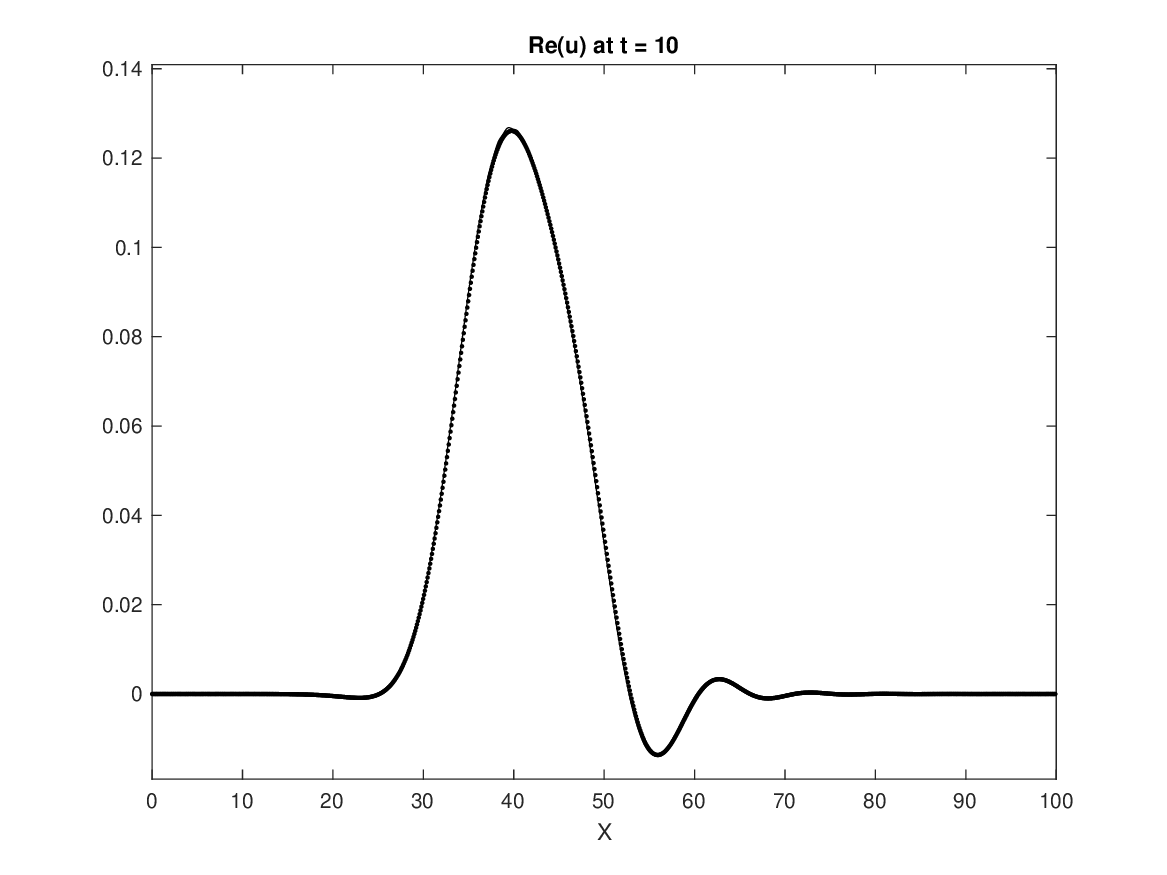}
\caption{Real part of the solution of the problem \eqref{schro1}-\eqref{periodic_cond} for $\alpha = 0$, $\gamma = -1$, $\eta = \beta =1$. Solid line: numerical solution. Pointed line: exact solution. }
\label{real_u_comp}
\end{figure}

\begin{figure}[h]
\includegraphics[width=15cm, height=10cm]{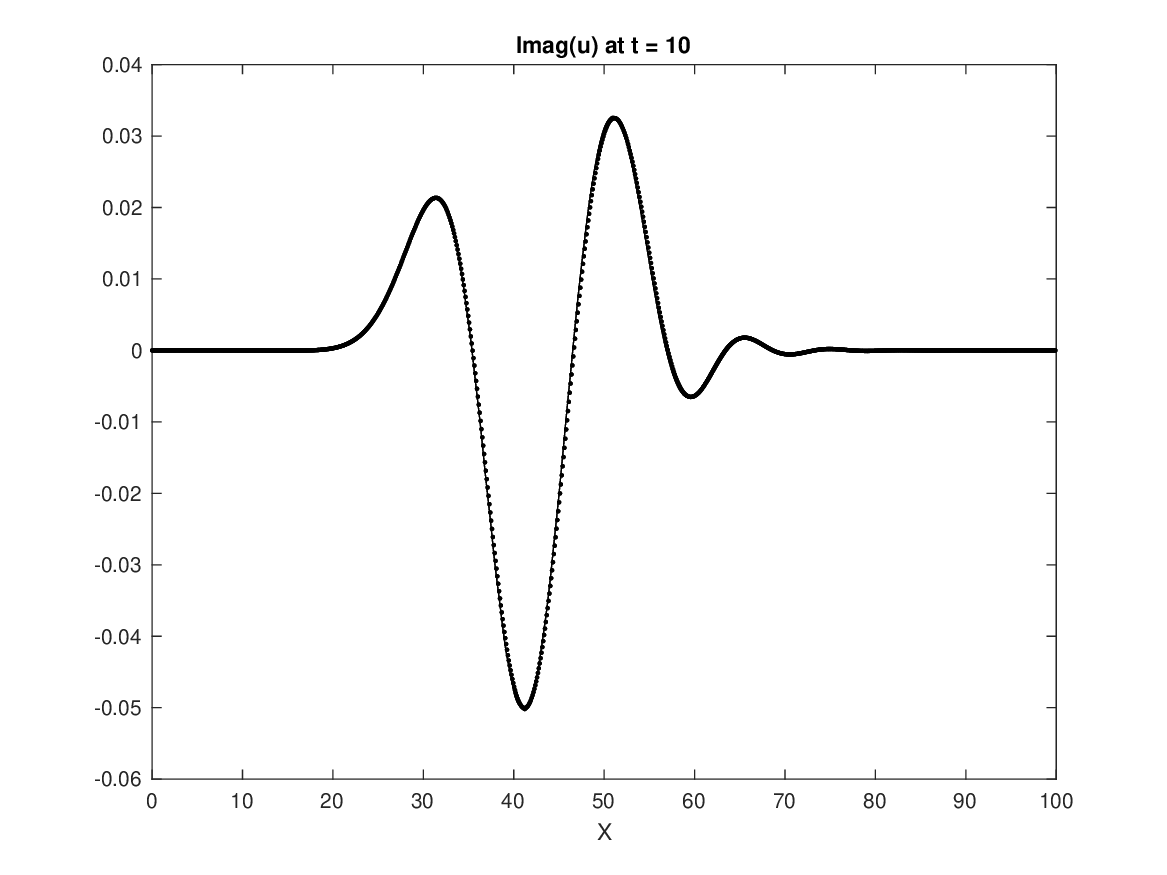}
\caption{Imaginary part of the solution of the problem \eqref{schro1}-\eqref{periodic_cond} for $\alpha = 0$, $\gamma = -1$, $\eta = \beta =1$. Solid line: numerical solution. Pointed line: exact solution. }
\label{imag_u_comp}
\end{figure}

As the first test, we validate the numerical solver using an exact solution of the linearized problem \eqref{schro1}-\eqref{periodic_cond}, i.e., $\alpha = 0$. In this case, as shown in the previous section, the exact solution is given by
\[
u(x,t) = \mathcal{F}^{-1} \Big( \hat{u}_{0,k}\text{exp}( -i w_k^2 - \eta w_k^2 - \beta i w_k^3 + i \gamma w_k )t \Big).
\]
We set the domain length as $L = 100$, the time step $\Delta t = 10/100$, the number of Fourier modes $N = 2^{10}$, 
and choose the parameters $\gamma = -1$, and $\eta = \beta =1$. The initial condition is taken to be a Gaussian pulse
\[
u(x,0) = u_0(x) = e^{-(x-30)^2}.
\]
Figures \ref{real_u_comp} and \ref{imag_u_comp} display the real and imaginary parts of the numerical solution superimposed with the exact solution. As shown, the numerical profiles closely match the exact ones, with high accuracy. No spurious dispersion or dissipation is observed in the simulations.

Next, we analyze the numerical error of the proposed scheme, as reported in Tables \ref{timeerror} and \ref{spaceerror}. The results confirm that the method achieves second-order convergence in time. In contrast, the spatial error of the Fourier spectral method decreases rapidly, exhibiting an algebraic convergence rate of approximately $N^{-4}$, where $N$ is the number of FFT points used in the computational domain.

In these experiments, the initial condition is
\[
u(x,0) = e^{-(x - 25)^2},
\]
and the model parameters are set as $\alpha = -1$, $\beta=0.5$, $\gamma = 0.5$, and $\eta = 0.5$.

\begin{table}[!h]
\begin{center}
\begin{tabular}{|c|c|c|c|} 
  \hline
  \multicolumn{4}{ |c| }{\textbf{Error in time}} \\ \hline   \hline
  \textbf{$\Delta t$} & \textbf{$\lVert U _{\Delta t}-U _{\Delta t/2}\rVert _{\infty}$} & \textbf{$\frac{\lVert U _{\Delta t}-U _{\Delta t/2}\rVert _{\infty}}{\lVert U _{\Delta t}\rVert _{\infty}}$} & \textbf{$\lambda_{1}$} \\
  \hline
$T/2$ & $2.5852 \times 10^{-5}$ & $5.2623 \times 10^{-}5$ &  2.0041 \\ \hline
$T/4$ &  $6.4447 \times 10^{-6}$ & $1.3119\times 10^{-5}$ & 2.0008 \\ \hline
$T/8$ &  $1.6102 \times 10^{-6}$ &  $3.2778\times 10^{-6}$ & 2.0002\\ \hline
$T/16$ &  $4.0251 \times 10^{-7}$ &  $8.1933 \times 10^{-7}$ & 2.0000 \\ \hline
$T/32$ &  $1.0063 \times 10^{-7}$ &  $2.0483 \times 10^{-7}$ & 2.0000 \\  \hline
\end{tabular} 
\end{center} 
\bigskip
\caption{Temporal convergence order of the Fourier-spectral method. Here $U_{\Delta t}$
denotes the approximation of the solution of equation \eqref{schro1} for the time step $\Delta t$, using the numerical solver \eqref{fullscheme}, with $N = 2^{20}$ FFT points fixed and refining progressively the time step. The final time is $T=1$. }
\label{timeerror} 
\end{table}

\begin{table}[!h]
\begin{center} 
	\begin{tabular}{|c|c|c|c|} 
		\hline
		\multicolumn{4}{ |c| }{\textbf{Error in space}} \\ \hline   \hline
		\textbf{$\Delta x$} & \textbf{$\lVert U _{\Delta x}-U _{\Delta x/2}\rVert _{\infty}$} & \textbf{$\frac{\lVert U _{\Delta x}-U _{\Delta x/2}\rVert _{\infty}}{\lVert U _{\Delta x}\rVert _{\infty}}$} & \textbf{$\lambda_{2}$} \\
		\hline
		$L/2^{8}$ & $    2.6727\times 10^{-6}$ & $5.4411 \times 10^{-6}$ &   4.0454 \\ \hline
		$L/2^{9}$ &  $  1.6187\times 10^{-7}$ & $3.2954\times 10^{-7}$ &  4.0122 \\ \hline
		$L/2^{10}$ &  $ 1.0032\times 10^{-8}$ &  $2.0423\times 10^{-8}$ &   4.0034\\ \hline
		$L/2^{11}$ &  $ 6.2550\times 10^{-10}$ &  $1.2732\times 10^{-9}$ & 4.0007\\ \hline
		$L/2^{12}$ &  $ 3.9075\times 10^{-11}$ &  $7.9539 \times 10^{-11}$ &   4.0002 \\ \hline
		$L/2^{13}$ & $2.4419\times 10^{-12}$& $4.9705\times 10^{-12}$& 3.9992 \\ \hline
		$L/2^{14}$ & $1.5270\times 10^{-13}$& $3.1082\times 10^{-13}$& 3.9793 \\ \hline
	\end{tabular} 
\end{center} 
\bigskip
\caption{Spatial convergence order of the Fourier-spectral method. Here $U_{\Delta x}$
denotes the approximation of the solution of equation \eqref{schro1} for the space size $\Delta x$, using the numerical solver \eqref{fullscheme} with $\Delta t = 1e-4$ fixed and refining progressively the spatial mesh size. The final time is $T=1$.}
\label{spaceerror}
\end{table}

We now turn to a numerical investigation of the behavior of solutions to the initial value problem \eqref{schro1}-\eqref{periodic_cond} in the limiting regime $\eta \to 0$. For this set of simulations, we fix the domain length $L = 50$, time step $\Delta t = 1.5/100$, and the number of Fourier modes $N = 2^9$. The parameters are set as $\gamma = 0$, $\alpha = -1$, and $\beta = 0$, and the initial condition is chosen as
We now analyze numerically the behavior of solutions of problem \eqref{schro1}-\eqref{periodic_cond} in the limit $\eta \to 0$. Here $L=50$, $\Delta t = 1.5/100$, $N=2^9$, $\gamma = 0$, $\alpha = -1$, $\beta = 0$, and the initial condition is
\[
u(x,0) = u_0(x) = e^{-(x - 25)^2}.
\]
In the nonlinear regime, i.e., when $\alpha \neq 0$, an exact analytical solution is not available. Therefore, we rely on the numerical solver developed in Section 4 to observe and interpret the qualitative behavior of the solutions in selected asymptotic parameter regimes.

\begin{figure}[h]
\includegraphics[width=15cm, height=15cm]{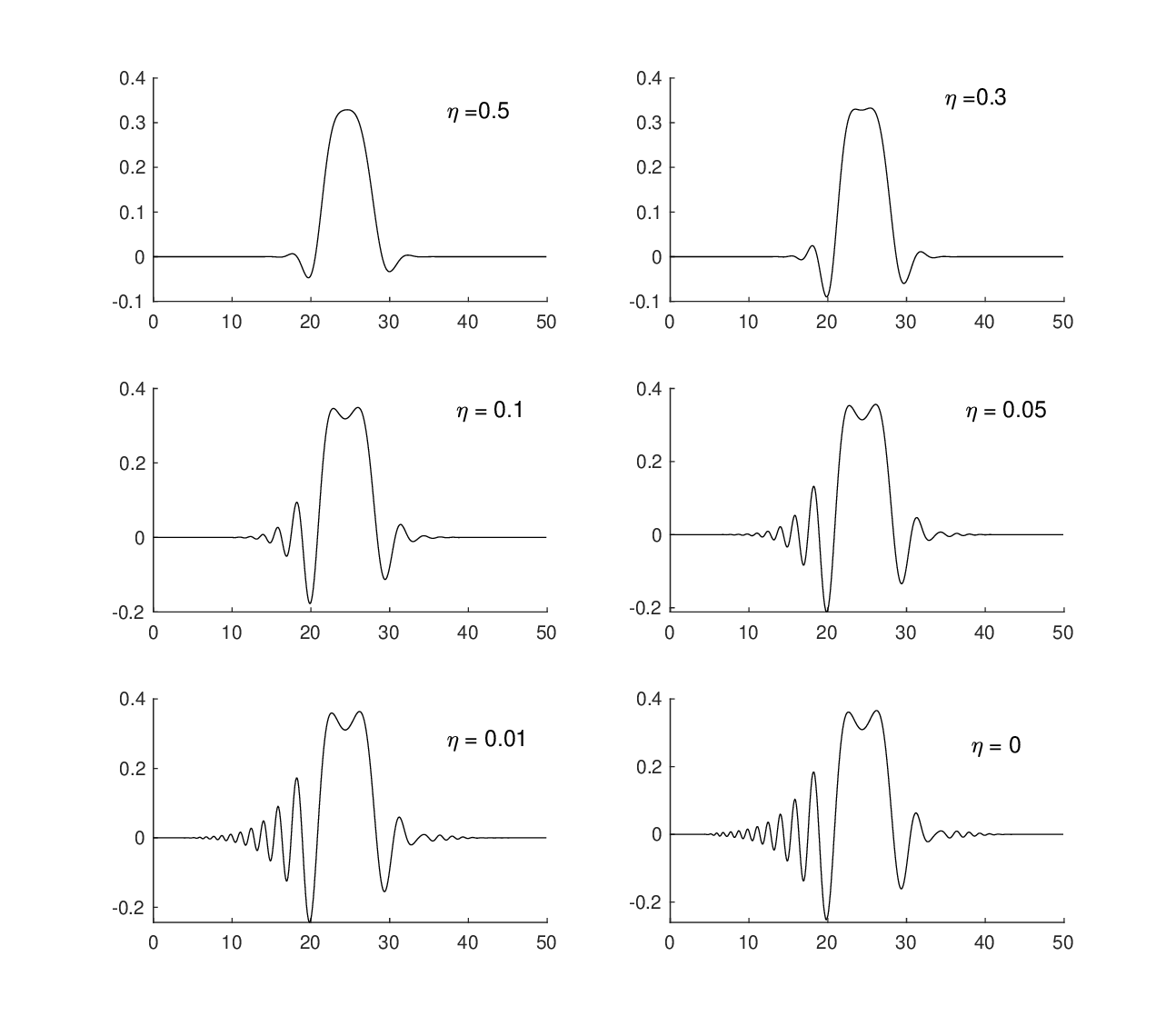}
\caption{The solution of problem \eqref{schro1}-\eqref{periodic_cond} with $\alpha =-1$, $\beta=0$, $\gamma = 0$, and a decreasing sequence of values of the parameter $\eta$.  }
\label{limit_eta}
\end{figure}

The influence of the diffusive term, governed by the parameter $\eta$, is illustrated in the sequence of plots shown in Figure \ref{limit_eta}. These simulations highlight the smoothing effect induced by diffusion: as $\eta$ increases, the oscillatory structures surrounding the main pulse diminish in amplitude and eventually vanish. This behavior reflects the dissipative nature of the diffusion term, which suppresses high-frequency components of the solution.

Conversely, as $\eta \to 0$, the solutions of the problem converge to the solution corresponding to the nondissipative case $\eta = 0$. This is in agreement with the theoretical results established in Section 3, where it was shown that the solution depends continuously on the parameter $\eta$ in appropriate function spaces. The numerical evidence further supports this theoretical prediction by demonstrating that the qualitative structure of the solution remains consistent and tends smoothly toward the zero-diffusion limit.

\begin{figure}[h]
\includegraphics[width=15cm, height=15cm]{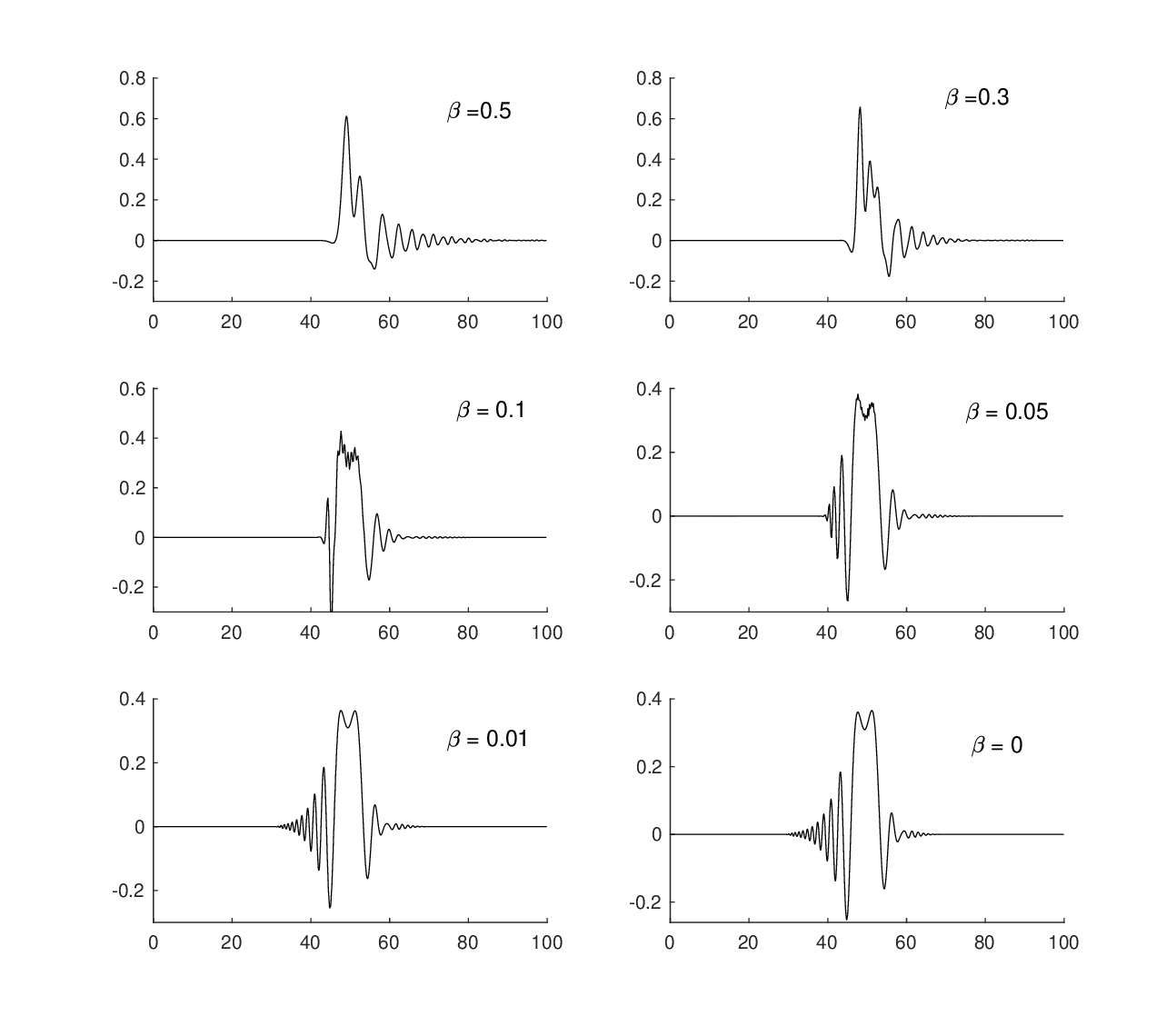}
\caption{The solution of problem \eqref{schro1}-\eqref{periodic_cond} with $\alpha =-1$, $\eta=0$, $\gamma = 0$, and a decreasing sequence of values of the parameter $\beta$.}
\label{limit_beta}
\end{figure}

We next consider the complementary asymptotic regime, in which the dispersive parameter $\beta$ tends to zero. For these experiments, we set $L = 100$, $N = 2^{10}$, and fix $\eta = 0$, while keeping the other parameters unchanged. The initial condition is given by
\[
u(x,0) = u_0(x) = e^{-(x-50)^2}.
\]
Figure \ref{limit_beta} displays the numerical solutions of problem \eqref{schro1}-\eqref{periodic_cond} for several values of the parameter $\beta$. As with the diffusion term, we observe that the solution structure evolves smoothly as $\beta$ decreases. For larger values of $\beta$, the solution exhibits dispersive features such as oscillatory tails; however, as $\beta \to 0$, these dispersive effects gradually vanish, and the solution approaches the profile corresponding to $\beta = 0$.

These numerical results confirm the theoretical analysis presented in Section 3, where it was shown that the solution of the model depends continuously on $\beta$. In particular, the simulations validate the asymptotic consistency of the numerical scheme and reinforce the understanding of the qualitative influence of the dispersive and diffusive terms in the governing equation.

\section{Conclusions}

In this work, we have analyzed a nonlinear Schr\"odinger-type equation that includes derivative nonlinearities and diffusion effects features that make the model relevant for several physical applications, such as low-order magnetization in ferromagnetic nanocables and the interaction of ferromagnetic solitons in weakly magnetized media.

We first established a local well-posedness result for the associated Cauchy problem, providing a rigorous mathematical foundation for the model. In addition, we proposed and implemented a Fourier spectral scheme tailored to the periodic setting, enabling accurate numerical approximations of the equation's solutions. Numerical experiments confirm the effectiveness of this method in capturing key dynamical behaviors.

We also explored the behavior of the model under various limiting regimes of its parameters, offering insights into the transition between different physical scenarios and aiding in the interpretation of solution structures.

Future work may include the study of global well-posedness, long-time dynamics, and the development of structure-preserving numerical methods. Further investigation into the physical implications of the model in more complex or higher-dimensional settings could also provide valuable extensions.


\end{document}